
\documentstyle{amsppt}
\topmatter
\author S.~Giller
\footnote[*]{supported by University of \L \' od\' z
grant No 458 \hfill{ }},
 C.~Gonera$^*$, P.~Kosi\' nski$^*$,
P.~Ma\' slanka
\footnote[**]{supported by KBN grant No 2P30221706p02 \hfill{ }},
\endauthor
\title The Quantum Galilei Group \endtitle
\address
S.~Giller, C.~Gonera, P.~Kosi\' nski
\newline Department of Field Theory
\newline University of \L \' od\' z
\newline Pomorska 149/153
\newline 90-236 \L \' od\' z, Poland
\newline
\newline
P.~Ma\' slanka
\newline Institute of Mathematics
\newline University of \L \' od\' z
\newline Banacha 22
\newline 90-238 \L \' od\' z, Poland
\endaddress
\abstract
The quantum Galilei group $G_{\varkappa}$ is defined.
The bicrossproduct structure of $G_{\varkappa}$ and the corresponding
Lie algebra is revealed. The projective representations for
twodimensional quantum Galilei group are constructed.
\endabstract
\endtopmatter
\NoRunningHeads
\magnification=\magstep1
\document

\TagsOnRight
\define\ti{\widetilde}
\define\ka{\varkappa}
\define\ep{\varepsilon}
\define\ro{\varrho}

\head I.~Introduction \endhead

Some attention has been recently paid to the so-called deformed
Poincar\' e algebra \cite{1} and group \cite{2}.
They are interesting because they provide relatively mild
deformation of classical space-time symmetries depending on
dimensionful deformation parameter. It is therefore interesting
to study in more detail some properties of the resulting structure.
In particular, in Ref.~\cite{3} some preliminary remarks have
been made concerning the non-relativistic limit of deformed symmetry.

In the present paper we continue this analysis. We construct the
$\varkappa$-deformed Galilei group and reveal the bicrossproduct
structure of both the algebra and the group (cf. Ref.~\cite{4}).

We address also to the question of the quantum counterpart of
projective representations. They are defined in analogy with
classical case. For twodimensional deformed Galilei group
the projective representations are explicitly constructed.

\head II.~The deformed Galilean algebra \endhead

The deformed Galilean algebra $\frak g_{\ka}$ can be obtained
by contraction procedure from $k$-Poincar\' e algebra in the way
described in Ref.~\cite{3}. We make the rescaling
$P_0\to P_0 \cdot c$, $L_i\to L_i \cdot c^{-1}$ and let
$c\to \infty$ keeping $kc\equiv \ka =const$. The resulting structure
reads~:
$$\aligned
[M_i, P_j]=i \ep_{ijk}P_k & ,\qquad [M_i,P_0]=0, \\
[M_i, M_j]=i \ep_{ijk}M_k & ,\qquad [P_{\mu},P_{\nu}]=0, \\
[M_i,L_j]=i \ep_{ijk} L_k & ,\\
[L_i,P_0]=i P_i & ,\qquad [L_i, P_k]=0, \\
[L_i,L_j]=\frac{1}{4\ka^2}\ep_{ijk} P_k
\left( \vec{P} \cdot \vec{M} \right)  & ,\\
\triangle M_i=M_i\otimes I+I\otimes M_i & ,\qquad
\triangle P_0=P_0\otimes I+I\otimes P_0, \\
\triangle P_i=P_i\otimes e^{-\frac{P_0}{2\ka}}
              +e^{\frac{P_0}{2\ka}}\otimes P_i & , \\
\triangle L_i=L_i\otimes e^{-\frac{P_0}{2\ka}}
              +e^{\frac{P_0}{2\ka}}\otimes L_i &
              -\frac{1}{2\ka}\ep_{ikl}
    \left( e^{\frac{P_0}{2\ka}}M_k\otimes P_l
    +P_k\otimes e^{-\frac{P_0}{2\ka}}M_l \right) , \\
S(P_{\mu})=-P_{\mu} & ,\qquad S(M_i)=-M_i, \\
S(L_i)=-L_i-\frac{3i}{2\ka}P_i & .
\endaligned \tag1 $$

Let us note that this algebra is obtained by contraction
accompanied with strong deformation limit $k\to 0$. It
seems that there exists no nonrelativistic limit
($c\to \infty$) with $k$ kept fixed \cite{5}.

The Casimir operators can be also obtained by contraction. They read~:
$$\aligned
& C_1=\vec{P}^2, \\
& C_2=\frac{\vec{P}^2}{4\ka^2} \left( \vec{P} \cdot \vec{M} \right)^2
+\left( \vec{P} \times \vec{L} \right)^2 .
\endaligned \tag2$$

Obviously, in the limit $\ka \to \infty$ we recover standard Galilean
structure. As in the case of $k$-deformed Poincar\' e algebra \cite{4}
one can show that our algebra has a bicrossproduct structure.
To this end we define new generators~:
$$\aligned
\ti{M}_i=M_i, & \qquad \ti{P}_0=P_0, \\
\ti{P}_i=P_i e^{-\frac{P_0}{2\ka}}, & \qquad
\ti{L}_i=L_i e^{-\frac{P_0}{2\ka}}
-\frac{\ep_{ijk}}{2\ka}M_jP_k e^{-\frac{P_0}{2\ka}}.
\endaligned \tag3$$
Then we obtain
$$\frak g_{\ka}=T \bowtie
U\left( \ti{M}, \ti{L} \right) ,\tag4$$
where
$U\left( \ti{M}, \ti{L} \right) $ is universal covering of Lie algebra
$\frak e(3)$
$$\aligned
[\ti{M}_i,\ti M_j]=i\ep_{ijk}\ti{M}_k, \qquad
[\ti{M}_i,\ti{L}_j]= & i \ep_{ijk} \ti{L}_k, \qquad
[\ti{L}_i,\ti{L}_j]=0, \\
\triangle \ti{M}_i=\ti{M}_i\otimes I+I\otimes \ti{M}_i, \qquad &
\triangle \ti{L}_i=\ti{L}_i\otimes I+I\otimes \ti{L}_i
\endaligned \tag5$$
while $T$ is defined by
$$\aligned
& [\ti{P}_{\mu},\ti{P}_{\nu}]=0, \\
& \triangle \ti{P}_0=\ti{P}_0\otimes I+I\otimes \ti{P}_0, \qquad
\triangle \ti{P}_i=\ti{P}_i\otimes e^{-\frac{\ti{P}_0}{\ka}}
+I\otimes \ti P_i.
\endaligned \tag6$$
Equation (4) is readily verified once we define
$$\aligned
& \ti{M}_i \triangleright \ti{P}_0=0, \qquad
\ti M_i\triangleright \ti P_j=i\ep_{ijk}\ti P_k, \qquad
\ti L_i\triangleright \ti P_0=i\ti P_i, \\
& \ti L_i\triangleright\ti P_j=\frac{i}{2\ka}\delta_{ij}
\vec{\ti P}^2-\frac{i}{\ka}\ti P_i\ti P_j, \\
& \delta \left( \ti M_i \right) =\ti M_i\otimes I, \qquad
\delta \left( \ti L_i \right) =\ti L_i \otimes
e^{-\frac{\ti P_0}{\ka}}-\frac{i}{\ka} \ep_{ijk} \ti M_j
\otimes \ti P_k.
\endaligned \tag7$$

Let us now define the cocommutator
$$\sigma=(\triangle -\tau \circ \triangle )
\left( \operatorname{mod} \frac{1}{\ka} \right) ; \tag8$$
we obtain
$$\aligned
& \sigma (P_0)=0, \qquad \sigma (M_i)=0, \\
& \sigma (P_i)=-\frac{1}{\ka}P_i\wedge P_0, \\
& \sigma (L_i)=-\frac{1}{\ka} (L_i\wedge P_0 +\ep_{ikl} M_k\wedge P_l ).
\endaligned \tag9$$

However, $\sigma$ is not a coboundary -- the classical $r$-matrix does
not exist. This can be verified by direct calculations.

\head III.~The deformed Galilei group \endhead

The deformed Galilei group $G_{\ka}$ can be obtained
either via contraction procedure or by quantizing the
Poisson structure on classical group implied by equations (9).
First we shall consider contraction procedure.

The relations defining classical Lorentz group,
$g_{\mu \nu}\Lambda_{\alpha}^{\mu} \Lambda_{\beta}^{\nu}
=g_{\alpha \beta}$, can be solved explicitly in terms of
rotation matrix $R_j^i$ and velocity $v^i$.
$$\aligned
& \Lambda_0^0=\left( 1+\frac{\vec{v}^2}{c^2} \right)^{\frac{1}{2}}, \\
& \Lambda_0^i=\frac{v^i}{c}, \\
& \Lambda_i^0=\frac{v^kR_i^k}{c}, \\
& \Lambda_i^k=\left( \delta_l^k+ \left( \left( 1+\frac{\vec{v}^2}{c^2}
\right)^{\frac{1}{2}} -1 \right) \frac{v^k v^l}{\vec{v}^2} \right)
R_i^l .
\endaligned \tag10$$
If we further redefine the translation sector, $a^i\to a^i$,
$a^0\to c\tau$, and let $c\to \infty$, $kc\equiv \ka =const$
in the commutation rules and coproduct defining the
$k$-Poincar\' e group \cite{2} we arrive at the following
structure
$$\aligned
& [R_j^i,R_l^k]=0, \qquad [R_j^i,v^k]=0, \qquad [v^i,v^k]=0, \\
& [a^i,a^j]=0, \qquad [\tau ,a^i]=\frac{i}{\ka}a^i, \\
& [\tau ,v^i]=\frac{i}{\ka}v^i, \qquad
  [v^i,a^j]=-\frac{i}{\ka} \left( v^iv^j
  -\frac{1}{2} \vec{v}^2 \delta_{ij} \right)  ,\\
& [R_j^i,\tau ]=0, \qquad [R_j^i, a^k]=-\frac{i}{\ka}
  \left( v^iR_j^k -\delta_{ik} v^mR_j^m \right) , \\
& \triangle R_j^i=R_k^i \otimes R_j^k, \\
& \triangle v^i=R_j^i \otimes v^j +v^i\otimes I , \\
& \triangle \tau =\tau \otimes I+I\otimes \tau ,\\
& \triangle a^i=R_j^i\otimes a^j+v^i\otimes \tau+a^i\otimes I .
\endaligned \tag11$$

The same structure is obtained by quantizing the Poisson structure on classical
Galilei group  implied by the cocommutator $\sigma$ (equation (9)). We define
the Poisson bracket by
$$\langle \{ f,g\},X\rangle =-i \langle f\otimes g, \sigma (X) \rangle ,
\tag12$$
where $\langle \ ,\ \rangle$ defines the classical duality between
Lie group and algebra which, in our case, reads
$$\aligned
\langle \tau ,P_0\rangle =i, \qquad
& \langle v^i, L_k\rangle =-i \delta_k^i, \\
\langle a^i,P_k \rangle =-i\delta_k^i, \qquad
& \langle R_j^i, M_k\rangle =-i\ep_{ijk} ,
\endaligned \tag13$$
the remaining brackets being zero. A straightforward calculations (see
Appendix)
allow us to recover the commutation rules (11).

Again, one can equipp $G_{\ka}$ with bicrossproduct structure.
Namely
$$G_{\ka}=T^*\bowtie C(E(3)), \tag14$$
where $C(E(3))$ is the algebra of functions on classical
group $E(3)$, generated by $R_j^i$ and $v^i$ while
$T^*$ is defined by
$$\aligned
[\tau ,a^i]=\frac{i}{\ka} a^i, \qquad
& [a^i,a^j]=0, \\
\triangle a^i=a^i\otimes I+I\otimes a^i, \qquad
& \triangle \tau =\tau \otimes I+I\otimes \tau ;
\endaligned \tag15$$
moreover,
$$\aligned
& f\triangleleft g=[f,g], \qquad f\in C(E(3)), \quad g\in T^*, \\
& \beta (\tau )=I\otimes \tau , \\
& \beta (a^i)=R_j^i\otimes a^j +v^i\otimes \tau .
\endaligned \tag16$$

It is likely that the recent proof of duality for $k$-Poincar\' e
algebra and group \cite{6}, \cite{7} can be adapted to our pair
$(G_{\ka}, \frak g_{\ka})$.

As next step we define the Galilean space-time. It is algebra
generated by four elements $t$, $x^i$, subject to the following
commutation rules
$$[t,x^i]=\frac{i}{\ka}x^i, \qquad [x^i, x^j]=0. \tag17$$

The Galilei group $G_{\ka}$ acts covariantly
according to the rules
$$\aligned
t\to & I\otimes t+\tau\otimes I, \\
x^i\to & R_j^i\otimes x^j +v^i\otimes t +a^i\otimes I .
\endaligned \tag18$$

\head IV.~Projective representations \endhead

It is well known that the Galilei group posses one-parameter
family of projective representations \cite{8}. Moreover,
they are exactly those representations which have physical meaning;
the parameter labelling inequivalent projective representations
is simply the mass of a particle. So the natural question arises
whether the projective representations have
some counterpart in the deformed case.

We have at our disposal no general theory of projective
representations for quantum groups. However, we can try to follow closely
the classical case and define an unitary projective represention of the
quantum group $\Cal A$ as a map $\ro :H\to H\otimes \Cal A$
satisfying
$$(\ro \otimes I)\circ \ro (\psi)=(I\otimes \omega )
((I\otimes \triangle)\circ \ro (\psi )), \qquad \psi \in H ;
\tag19$$
here $H$ is the relevant Hilbert space of states and $\omega$
is an unitary element of $\Cal A\otimes \Cal A$. As in the classical case
one easily derives the following consistency condition on $\omega$
(implied by (co-)associativity)
$$(\omega \otimes I)(\triangle \otimes I)\omega
=(I\otimes \omega )(I\otimes \triangle ) \omega .\tag20$$

Two projective representations $\ro$, $\ti \ro$, are called
equivalent if there exists an unitary element $a\in \Cal A$
such that
$$\ti \ro =(I\otimes a)\ro .\tag21$$

The corresponding multipliers  are related by the formula
$$(a\otimes a)\omega =\ti \omega \triangle (a). \tag22$$

Two such multipliers will be called equivalent. Obviously,
a multiplier $\omega$ is trivial (the representation is
equivalent to the vector one) if
$$(a\otimes a)\omega =\triangle (a) .\tag23$$

So, the problem reduces to that of finding solutions
of equation (20) which are not of the form (23). We shall solve
this problem for twodimensional Galilei group as a toy model.

The dimensional reduction applied to $G_{\ka}$ gives rise
to the following structure
$$\aligned
& [\tau ,a]=\frac{i}{\ka}a, \qquad [\tau ,v]=\frac{i}{\ka}v,
  \qquad [v,a]=-\frac{i}{2\ka}v^2, \\
& \triangle v=v\otimes I+I\otimes v, \qquad
  \triangle \tau =\tau \otimes I+I\otimes \tau, \\
& \triangle a=a\otimes I+I\otimes a+v\otimes \tau .
\endaligned \tag24$$

For the the classical group the multiplier takes the form
$$\omega_0=e^{i\varphi_0}, \qquad
\varphi_0=-m\left( \frac{v^2}{2}\otimes \tau +v\otimes a \right)
\tag25$$

Therefore we assume the same exponentional form for our quantum
multiplier
$$\omega =e^{i\varphi} \tag26a$$
and expand $\varphi$ in inverse powers of $\ka$ :
$$\varphi =\sum_{n=0}^{\infty} \frac{\varphi_n}{\ka^n} \tag26b$$

To provide nontriviality of $\omega$ we choose the first term in the
expansion to be given by equation (25).

Inserting (26a), (26b) into equation (20) and comparying the coefficients
in front of $\dfrac{1}{\ka}$ we arrive at the equation for $\varphi_1$.
The contribution to both sides are twofold~:
\roster
\item"{(i)}" the ones coming from $\varphi_1$ and
\item"{(ii)}" those coming from the commutators of $\varphi_0$
              due to Hausdorff formula.
\endroster
Let us note that after taking commutators all elements can be viewed
as classical ones -- the contributions due to the noncommutativity
are of higher orders in $\dfrac{1}{\ka}$. The resulting
equation reads
$$\multline
\varphi_1\otimes I-I\otimes \varphi_1 +(\triangle \otimes I)\varphi_1
-(I\otimes \triangle )\varphi_1 \\
=\frac{i\ka}{2}([I\otimes \varphi_0 ,(I\otimes \triangle )\varphi_0]
-[\varphi_0 \otimes I,(\triangle \otimes I) \varphi_0 ]) ,
\endmultline \tag27$$
where, after calculating the commutators on right-hand side
the equation can viewed as defined on classical Galilei group.

With some effort we can solve  equation (27) (which is typical
cohomological equation) for $\varphi_1$. The particular solution reads
$$\varphi_1 =\left( -\frac{1}{4}\frac{mv^2}{\ka} \otimes I \right) \varphi_0 .
\tag28$$

Therefore we assume the following Ansatz for $\varphi$~:
$$\varphi =\left( f\left( v^2 \right) \otimes I \right) \varphi_0 .
\tag29$$

By inserting (29) and (26a) into equation (20) we find
after some tedious calculations that this Ansatz is
consistent and provides the following expression for $\omega$~:
$$\omega =
e^{-i \left( \frac{2\ka}{v^2} \ln \left( 1+\frac{mv^2}{2\ka} \right)
      \otimes I \right)
      \left( \frac{v^2}{2}\otimes \tau +v\otimes a \right) }
\tag30$$
or, equivalently
$$\omega =
e^{-i \ka \ln \left( 1+\frac{mv^2}{2\ka} \right) \otimes \tau }
e^{\frac{-imv}{1+\frac{mv^2}{2\ka}}\otimes a }. \tag31$$

In order to obtain the relevant representation we can argue as follows.
In the classical case, given the multiplier $\omega (g,g')$ one can
construct the relevant representation by considering
the linear space of functions defined over the group and defining
the group action by the formula : $f(g)\to \omega (g,g')f(gg')$.
Equation (20) allow us to conclude  that this is projective
representation, $\omega (g,g')$ being the relevant multiplier.
If $\omega (g,g')$, as a function of the first variable can
be viewed as defined over some coset space we can take $f$ to
be also defined over this coset space. This is the case for our
$\omega_0$; as a function of first variable it is defined over
the coset space parametrized by boosts.
If we call $v=-\frac{P}{m}$ we recover the standard representation.
The whole procedure can be applied to the quantum case; $\omega$ has
the same property as $\omega_0$ -- it depends on the first variable
only through boost $v$. In this way we arrive at the following form of
(unitary) representation for twodimensional quantum Galilei
group (24)~:
$$\ro :f(p)\to
e^{-i \ka \ln \left( 1+\frac{p^2}{2m\ka} \right) \otimes \tau }
e^{i\frac{p}{1+\frac{p^2}{2m\ka}}\otimes a }
f(p\otimes I-I\otimes mv) .\tag32$$

Let us note that this representation is well defined only
for $\ka >0$; the same phenomenon occurs for deformed
Poincare group \cite{9}.

The above procedure can be extended to fourdimensional case.
Moreover, using duality relations group$\leftrightarrow$algebra,
one can find the infitesimal form of representation and, consequently,
the algebra obeyed by generators.
This should indicate the way of constructing the
{\it quantum-mechanical} extension of quantum
Galilei group -- in analogy with classical case.
Another problem is how to obtain the representation
of quantum Galilei group from those of deformed Poincar\' e
group. All these questions will be addressed to in subsequent
publication.

\head Appendix \endhead

We shall show that the commutation rules (11) can be obtained
by quantizing the Poisson structure implied by commutator $\sigma$.

The classical duality group$\leftrightarrow$algebra
is defined as
$$\langle \Phi ,X\rangle=-i\frac{d}{dt}\Phi
\left( e^{itX} \right) |_{t=0} \tag{A1}$$
which implies equations (13).

We fix the following order of factors in monomials belonging
to the universal covering of Galilei algebra~: $\vec{M}$,
$\vec{L}$, $\vec{P}$, $P_0$. The standard duality rules
$$\langle \Phi , XY\rangle
=\langle \triangle \Phi ,X\otimes Y\rangle, \qquad
\langle \Phi \Psi ,X\rangle
=\langle \Phi \otimes \Psi, \triangle X\rangle$$
give, in particular
$$\aligned
& \langle R_{ij},J_{n_1}\ldots J_{n_k}X \rangle
  =(-i)^k\delta_{XI}\ep_{il_1n_1}\ep_{l_1l_2n_2}\ldots \ep_{l_{k-1}jn_k} ,\\
& \langle \tau , XH^k\rangle =i\delta_{XI}\delta_{k,1} ,\\
& \langle v_{i},J_{n_1}\ldots J_{n_k}L_m \rangle
  =-i \langle R_{im},J_{n_1}\ldots J_{n_k} \rangle , \\
& \langle a_{i},J_{n_1}\ldots J_{n_k}P_m \rangle
  =-i \langle R_{im},J_{n_1}\ldots J_{n_k} \rangle , \\
& \langle a_{i},J_{n_1}\ldots J_{n_k}L_mP_0 \rangle
  =\langle R_{im},J_{n_1}\ldots J_{n_k} \rangle .
\endaligned \tag{A2}$$

As an example let us calculate the Poisson bracket $\{ R_n^m,a^r\}$~:
$$\langle \{ R_n^m,a^r\}, X\rangle
=(-i)\langle R_n^m\otimes a^r, \sigma (X)\rangle . \tag{A3}$$

Using duality rules and the definition of cocommutator
we readily infer that the right-hand side of (A3) does not
vanish only for $X=M^kL_i$, where $M^k$ stands for the product
 of $k$ $M$'s with arbitrary indices. A straightforward
calculation based on equations (A1), (A2), (9) and the duality
rules gives
$$\langle \{ R_n^m,a^r\}, M^kL_i \rangle
=\frac{1}{\ka} \langle \delta_{mr} v^p R_n^p -v^mR_n^r,
M^kL_i\rangle .\tag{A4}$$

Now, we easily check that
$$\langle \delta_{mr}v^pR_n^p -v^mR_n^r, X\rangle =0$$
for $X\neq M^kL_i$; therefore
$$\{ R_n^m, a^r\}=-\frac{1}{\ka}
\left( v^mR_n^r-\delta_{mr}v^pR_n^p \right) .$$

\enddocument